\documentclass[a4paper,11pt]{article}
\usepackage{titlesec}
\usepackage{amsthm}
\usepackage{amssymb}
\usepackage{color}
\usepackage{mathtools}
\usepackage{geometry}
\usepackage{algorithm2e}

\usepackage{verbatim}
\numberwithin{equation}{section}

\newtheorem{theorem}{Theorem}[section]

\theoremstyle{definition}

\newtheorem{remark}[theorem]{Remark}
\newcommand{\bfs}[1]{{\boldsymbol #1}}

\usepackage{array}
\newcolumntype{C}[1]{>{\centering\arraybackslash}m{#1}}
\usepackage{graphics}
\usepackage{amssymb, latexsym, amsmath, amsfonts, epsfig}
\usepackage{bm}
\usepackage{graphicx}
\usepackage{color}
\usepackage{epstopdf}
\usepackage{comment}
\usepackage{soul}
\usepackage[normalem]{ulem}
\usepackage{float}
\usepackage{MnSymbol}

\usepackage{caption}
\usepackage{amsmath,amssymb}

\usepackage{lipsum}
\usepackage{amsfonts}
\usepackage{graphicx}
\usepackage{epstopdf}
\usepackage{algorithmic}
\ifpdf
  \DeclareGraphicsExtensions{.eps,.pdf,.png,.jpg}
\else
  \DeclareGraphicsExtensions{.eps}
\fi

\usepackage{enumitem}
\setlist[enumerate]{leftmargin=.5in}
\setlist[itemize]{leftmargin=.5in}

\title{A boundary-penalized isogeometric analysis for second-order hyperbolic equations}

\author{Quanling Deng\thanks{Corresponding author. School of Computing, Australian National University, Canberra, ACT 2601, Australia (Quanling.Deng@anu.edu.au) }
  \and Pouria Behnoudfar\thanks{Mineral Resources, Commonwealth Scientific and Industrial Research Organisation (CSIRO), Kensington, Perth, WA 6152, Australia}
  \and Victor Calo\thanks{School of Electrical Engineering, Computing and Mathematical Sciences, Curtin University, Perth, WA 6102, Australia} }

\usepackage{amsopn}

\begin{document}

\maketitle

% REQUIRED
\begin{abstract}
  Explicit time-marching schemes are popular for solving time-dependent partial differential equations; one of the biggest challenges these methods suffer is increasing the critical time-marching step size that guarantees numerical stability.  In general, there are two ways to increase the critical step size.  One is to reduce the stiffness of the spatially discretized system, while the other is to design time-marching schemes with larger stability regions.  In this paper, we focus on the recently proposed explicit generalized-$\alpha$ method for second-order hyperbolic equations and increase the critical step size by reducing the stiffness of the isogeometric-discretized system.  In particular, we apply boundary penalization to lessen the system's stiffness.  For $p$-th order $C^{p-1}$ isogeometric elements, we show numerically that the critical step size increases by a factor of $\sqrt{\frac{p^2-3p+6}{4}}$, which indicates the advantages of using the proposed method, especially for high-order elements.  Various examples in one, two, and three dimensions validate the performance of the proposed technique.
% \textbf{Mathematics Subjects Classification}:   65P99, 65M99, 76M28
\end{abstract}

% REQUIRED

% REQUIRED

\section{Introduction} \label{sec:intr} 

We consider the following second-order hyperbolic equation
\begin{equation}\label{eq:pde}
  \begin{cases}
    \begin{aligned}
      \ddot u(\bfs{x}, t) -  \nabla \cdot ( \kappa (\bfs{x}) \nabla
      u(\bfs{x}, t) )
      & =f(\bfs{x}, t),  &&  \bfs{x} \in \Omega,  \ t\in (0, T], \\
      u(\bfs{x}, t) & = u_D,  &&  \bfs{x} \in \partial\Omega, \ t>0, \\
      u(\bfs{x}, 0) & = u_0,  &&  \bfs{x} \in \Omega, \\
      \dot u(\bfs{x}, 0) & = v_0,  &&  \bfs{x} \in \Omega, \\
    \end{aligned}
  \end{cases}
\end{equation}
where $\Omega = [0,1]^d \subset \mathbb{R}^d, d=1,2,3$, is a bounded open domain with Lipschitz boundary $\partial\Omega$.  $\nabla \cdot$ is the divergence operator, $\nabla$ is the gradient operator, $0<\kappa_0 \le \kappa (\bfs{x}) \le \kappa_1$ is a diffusion coefficient representing the material property, $f$ is a forcing function, $u$ is the unknown, $u_D$ provides boundary data, $u_0$ and $v_0$ provide initial data. The superposed dots refer to the time-derivatives where $\dot u = \frac{\partial u}{\partial t}, \ddot u = \frac{\partial^2 u}{\partial t^2}$ denote the velocity and acceleration, respectively.  Equation~\eqref{eq:pde} arises in mathematical modeling of various engineering and scientific problems.  For example, in structural engineering, we model the structural vibrations due to dynamic excitations from wind, earthquakes, blasts, vehicular traffic, and operating machinery, as a system second-order hyperbolic differential equations~\cite{ courant2008methods, tedesco2000structural, craig2006fundamentals}.  In particular, undamped structural dynamics are usually modeled as~\eqref{eq:pde}; see, for example,~\cite{ paz2012structural}.  Equation~\eqref{eq:pde} is also referred to as the wave equation~\cite{ evans2010partial, serway2018physics, courant2008methods}, often seen in models where the underlying medium is homogeneous (when $\kappa$ is constant)~\cite{ serway2018physics}.  While the regularity (existence, uniqueness, and smoothness) of the solution of~\eqref{eq:pde} has been well-understood, finding its analytical solutions is impossible.  Therefore, the modeling equation~\eqref{eq:pde} with a general domain is usually solved numerically.

Various numerical methods have been developed to solve equation~\eqref{eq:pde}.  For the spatial discretization, isogeometric analysis (IGA)~\cite{ hughes2005isogeometric, cottrell2009isogeometric} is a state-of-the-art method that combines classical finite element analysis with computer-aided design and analysis tools.  A series of research works on IGA studies the advantages of the method over the classical finite element method (FEM) for spectral approximations~\cite{ nguyen2015isogeometric, cottrell2006isogeometric, hughes2008duality, hughes2014finite, calo2019dispersion, puzyrev2017dispersion, deng2018ddm, Garcia:2019, Hashemian:2021, deng2018dispersion} as well as to expand its applications to various engineering and scientific problems~\cite{ zhang2007patient, gomez2008isogeometric, buffa2010isogeometric, verhoosel2011isogeometric, wang2018structural, Sarmiento:2017, deng2019optimal, deng2018isogeometric, deng2022isogeometric, Dalcin:2016}.  We herein use IGA for the spatial discretization of problem~\eqref{eq:pde}.

\begin{figure}[h!]
\centering
\includegraphics[height=6cm]{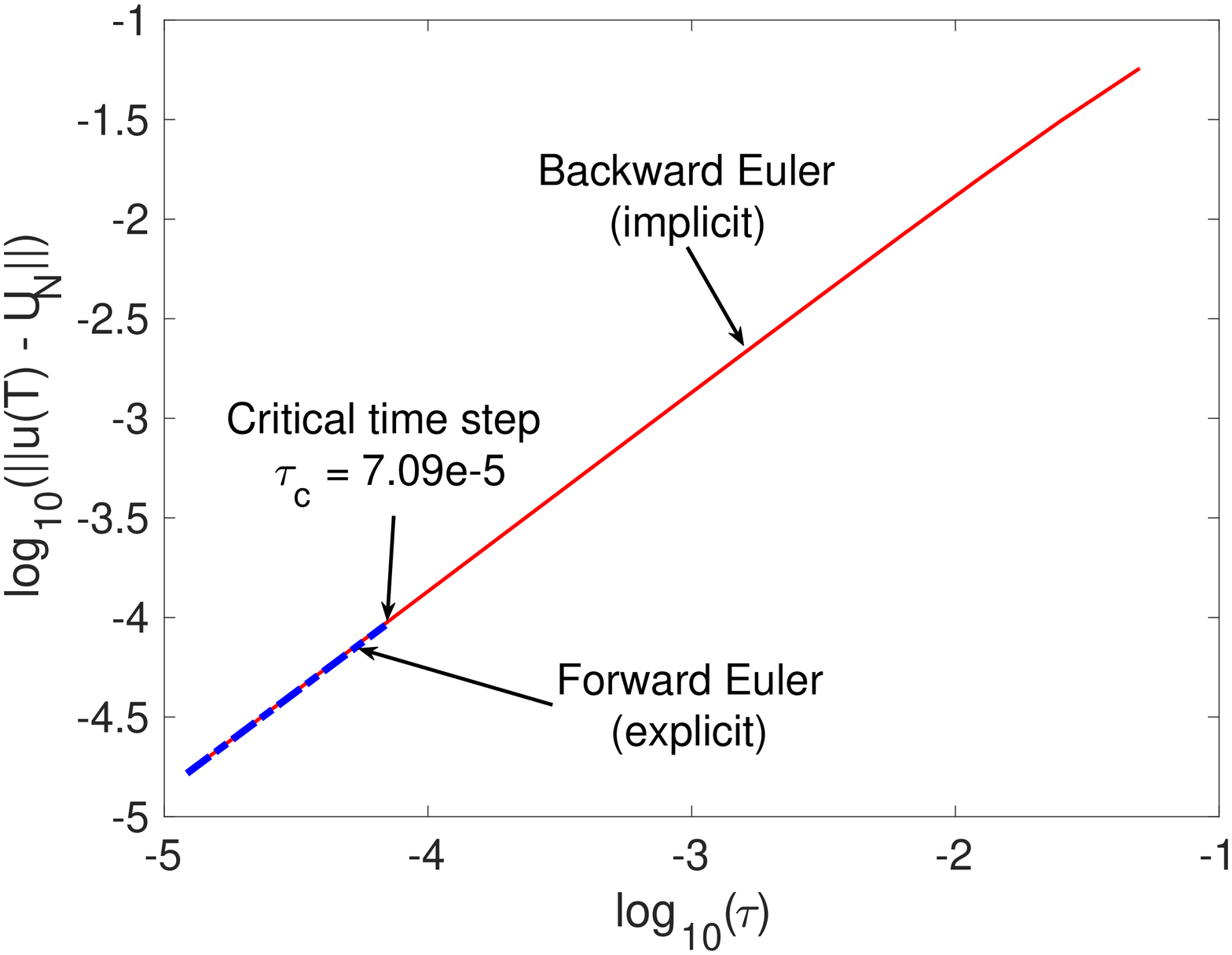} 
\caption{An example: $L^2$-norm errors versus time step sizes for both backward (implicit) and forward (explicit) Euler methods. The critical time step $\tau_c$ is the maximal time step size that ensures stability of the explicit scheme.}
\label{fig:ctime}
\end{figure}

There are explicit and implicit time-marching schemes to advance the problem in time~\cite{ butcher2016numerical}. The forward and backward Euler schemes are the simplest explicit and implicit schemes, respectively. For the second-order hyperbolic problem, especially in the field of structural dynamics, the generalized-$\alpha$ method~\cite{ chung1993time} is a widely-used implicit scheme. In this paper, we adopt the recently-developed explicit generalized-$\alpha$ methods~\cite{ behnoudfar2022explicit, labanda2021explicit} for the temporal discretization of problem~\eqref{eq:pde}.  While maintaining the second-order accuracy and other features of the original generalized-$\alpha$ method, this scheme also has the advantages of an implicit scheme. Generally speaking, implicit schemes are popular due to their unconditional stability. In addition, the numerical errors vary continuously with the time step size; see the red line of Figure~\ref{fig:ctime} as an example. Their main drawback is that one needs to invert large sparse algebraic systems, which can be computationally costly. Moreover, implementing an implicit scheme is generally harder than an explicit one. Explicit schemes are widely-used for the sake of simplicity and computational efficiency. In particular, they do not require inverting matrices. In the region where the time step size $\tau \le \tau_c$, explicit schemes are significantly cheaper than implicit ones to reach the same approximation accuracy. Figure~\ref{fig:ctime} shows an example of the critical time step size $\tau_c$ for the explicit forward Euler scheme. The scheme fails for larger time-step sizes, $\tau > \tau_c$ due to the lack of stability. Thus, less accurate approximations are not possible when using an explicit scheme by enlarging the time step size.

In this paper, we enlarge the critical time step size by reducing the high-frequency unphysical stiffness of the spatial discretization. For any explicit time-marching scheme, the stability condition boils down to a condition $\tau^2 \lambda \le C$ where $C>0$ is a constant. $\lambda \in \mathbb{C}$ is the largest eigenvalue in magnitude, which characterizes the stiffness of the spatial discretization. The critical time step size is then $\tau_c = \sqrt{C/\lambda}$, implying that a reduced stiffness (smaller $\lambda$) enlarges the critical time step size $\tau_c$. Our recent work~\cite{ deng2021boundary, deng2021outlier} introduces a boundary penalization technique that significantly reduces the stiffness of the IGA discretized system by eliminating the outliers. These outliers correspond to stopping bands in the finite element discrete spectrum. In the spectral approximation of second-order elliptic operators, this technique removes the outliers that were first observed in \cite{ cottrell2006isogeometric}. Recently, alternative outlier removal techniques~\cite{ hiemstra2021removal, manni2022application} build extra boundary conditions in the approximation space. The elimination of the outliers reduces the unphysical stiffness of the discretized system. We herein focus on the use of the boundary penalization technique developed in~\cite{ deng2021boundary}. We show numerically that the critical time step size increases by a factor of $\sqrt{\frac{p^2-3p+6}{4}}$ when using $p$-th order IGA elements for the wave equation (i.e., Equation~\eqref{eq:pde} with $\kappa=1$).

The rest of this paper is organized as follows.  Section~\ref{sec:idea} presents a numerical method for solving the problem~\eqref{eq:pde}; we start with an IGA spatial discretization and then introduce the explicit generalized-$\alpha$ method for time marching.  Section~\ref{sec:bp} describes the boundary penalization technique and discusses the advances in stiffness reduction and critical time step size increment.  Section~\ref{sec:num} collects numerical results that demonstrate the performance of the proposed method.  We focus on the study of the critical time step size increment.  Concluding remarks are presented in Section~\ref{sec:conclusion}.

\section{Fully discrete solver} \label{sec:idea}

First, we discretize~\eqref{eq:pde} in space using isogeometric analysis and then discretize the resulting ordinary differential matrix system using the explicit generalized-$\alpha$ method.

\subsection{Isogeometric discretization}

Before discretizing in space the problem~\eqref{eq:pde}, we obtain the weak form of the problem; from it we obtain the Galerkin formulation.  Multiplying the equation~\eqref{eq:pde} by a test function and integrating by parts, we obtain the weak form of the problem~\eqref{eq:pde}~\cite{ larsson2003partial}: for $ f \in L^2$, find $u(\cdot, t)\in H^1_0(\Omega)$ such that
\begin{equation} \label{eq:gf}
  a(\ddot u, w) + b(u,w) = \ell(w), \qquad \forall w \in H^1_0(\Omega),
\end{equation}
 where 
\begin{equation*}
  \begin{aligned}
    a(w, v) = (w, v), \qquad b(w, v) = (\nabla w, \nabla v), \qquad \ell(w) = (f, w), \quad \forall w, v \in H^1_0(\Omega).
  \end{aligned}
\end{equation*}
The method imposes the boundary conditions through the approximation space $H^1_0(\Omega)$ where the time-marching method enforces the initial conditions.  We introduce a finite-dimensional approximation to $H^1_0(\Omega)$ within the Galerkin framework; we adopt the isogeometric analysis for spatial discretization.

For simplicity, at the discrete level, we first partition the unit domain $\Omega=[0,1]^d, d=1,2,3,$ into a mesh with uniform tensor-product elements.  Let $E$ be a generic element and denote its collection as $\mathcal{E}_h$ such that $\overline \Omega = \cup_{E\in \mathcal{E}_h} E$.  We denote $h = \max_{E \in \mathcal{E}_h} \text{diameter}(E)$.   Isogeometric analysis uses the B-splines as basis functions in the Galerkin framework. In 1D, they are given by the Cox-de Boor recursion formula~\cite{ de1978practical, piegl2012nurbs} while in multiple dimensions they are given as tensor-products of the 1D functions.  We refer to~\cite{ buffa2010isogeometric, evans2013isogeometric, Sarmiento:2017, Dalcin:2016} for detailed constructions and define the approximation space as
\begin{equation} \label{eq:bs}
  V^h_p = \{ w \in C^{p-1}(\Omega): w|_{\partial\Omega} = 0, w|_E \in \mathbb{P}^p,  \forall E \in \mathcal{E}_h \} \subset H_0^1(\Omega),
\end{equation}
where $\mathbb{P}^p$ is the space of $p$-th order polynomials.  

At the spatial semi-discrete level, the isogeometric analysis of~\eqref{eq:pde} is to find $u^h(\cdot, t) \in V^h_p$ such that 
\begin{equation} \label{eq:vfh}
  a(\ddot u^h, w^h) + b(u^h, w^h) = \ell(w^h), \quad \forall \ w^h \in V^h_p,
\end{equation} 
which leads to the matrix problem
\begin{equation} \label{eq:mpf}
  M \ddot U + K U = F,
\end{equation}
where $F_k =\ell(\phi_p^k)$, $M_{kl} = a(\phi_p^k, \phi_p^l), K_{kl} = b(\phi_p^k, \phi_p^l),$ and $U$ is the vector of the coefficients of the basis functions $\phi_p^k$ in a solution representation. The matrix problem~\eqref{eq:mpf} is usually referred to as a semi-discretized system of the problem~\eqref{eq:pde} and it is a system of ordinary differential equations (ODEs).

\subsection{Explicit generalized-$\alpha$ method} \label{sec:ga}

For the temporal discretization, we adopt the explicit generalized-$\alpha$ method to solve the ODEs~\eqref{eq:mpf}.  We partition of the time interval $[0,T]$ with a grid $0 = t_0 < t_1 < \cdots < t_N = T$. The time-marching step size $\tau_n = t_n - t_{n-1}$.  Again, for simplicity, we assume a uniform step size $\tau$.  We approximate  $U(t_n), \dot{U}(t_n), \ddot{U}(t_n)$ as $U_n, V_n, A_n$, respectively.

Given $U_n, V_n, A_n$, the explicit generalized-$\alpha$ time-marching method is to find $U_{n+1}$, $V_{n+1}$, $A_{n+1}$ such that
\begin{equation} \label{eq:ga}
  \begin{aligned}
    M A_{n+\alpha_m} + K U_n  &= F_{n+\alpha_f}, \\
    V_{n+1} & = V_n + \tau_n A_n + \tau_n^2 \gamma \lsem A_n \rsem, \\
    U_{n+1} & = U_n + \tau_n V_n + \frac{\tau_n^2}{2} A_n + \tau_n^2 \beta \lsem A_n \rsem, \\
  \end{aligned}
\end{equation}
where
\begin{equation} \label{eq:ga0}
  \begin{aligned}
    F_{n+\alpha_f} & = F(t_{n+\alpha_f}) = F(t_n + \alpha_f \tau_n), \\
    \lsem A_n \rsem& = A_{n+1} - A_n, \\
    A_{n+\alpha_m} & = A_n + \alpha_m \lsem A_n \rsem
  \end{aligned}
\end{equation}
with initial conditions
\begin{equation} \label{eq:ga1}
  \begin{aligned}
    U_0 & = U(0), \\
    V_0 & = V(0), \\
    A_0 & = M^{-1} (F_0 - K U_0).
  \end{aligned}
\end{equation}
We use the parameter settings of~\cite{ labanda2021explicit}, which allow dissipation control and second-order accuracy.  These parameters are:
\begin{equation} \label{eq:rho}
  \begin{aligned}
    \alpha_f  = 0, \quad 
    \gamma = \frac12 + \alpha_m, \quad 
    \alpha_m  = \frac{2 - \rho}{\rho+1}, \quad 
    \beta  = \frac{3\rho - 5}{(\rho - 2)( \rho + 1 )^2},
  \end{aligned}
\end{equation}
where $0 \le \rho \le 1$ is a user-specified parameter that controls high-frequency dissipation.  When $\rho=1$, the scheme has the largest stability region and is equivalent to the second-order explicit central finite difference method~\cite{ wanner1996solving}. This formulation is also a special case of the Newmark family~\cite{ newmark1959method}. In this case, the conditional stability condition is as follows~\cite[Chapter~9]{ hughes2012finite}
\begin{equation}
  \tau \omega^h \le 2,
\end{equation}
where $\omega^h$ is the discrete frequency such that 
\begin{equation} \label{eq:mevp}
  KU = (\omega^h)^2 MU.
\end{equation}
The critical time step size is $\tau_c = 2/\omega^h_{\text{max}}$, which is inversely proportional to the maximal discrete frequency $\omega^h_{\text{max}}$.  In isogeometric discretizations of~\eqref{eq:pde}, the resulting matrix problem~\eqref{eq:mevp} has ``outliers" in the highest-frequency region.  These outliers are unphysical approximations to the maximal physical frequency $\omega_{\text{max}}$ from above (meaning $\omega^h_{\text{max}} > \omega_{\text{max}}$) with large errors.  The removal of these outliers reduces the approximation errors on the maximal discrete frequency. Consequently, their removal increases the critical time step size $\tau_c$, enlarging the stability regions for explicit time marching.

\section{Boundary penalization for isogeometric analysis} \label{sec:bp}

Following~\cite{ deng2021boundary}, we present a boundary penalization technique that removes outliers from the discrete frequency.  The method depends on the order of the isogeometric element, thus we first define
\begin{equation}
  \alpha = \lfloor \frac{p-1}{2} \rfloor =
  \begin{cases}
    \frac{p-1}{2}, & p \quad \text{is odd}, \\
    \frac{p-2}{2}, & p \quad \text{is even}.
  \end{cases}
\end{equation}
In one dimension, the boundary penalization technique for isogeometric analysis of~\eqref{eq:pde} discretizes~\eqref{eq:gf} to find $\tilde u^h(\cdot, t) \in V^h_p$ for all $f \in L^2(\Omega)$ such that
\begin{equation} \label{eq:vfhnew}
  \tilde a(\ddot{\tilde u}^h, w^h) + \tilde b(\tilde u^h, w^h) = \ell(w^h), \quad \forall \ w^h \in V^h_p,
\end{equation} 
where for $w, v \in V^h_p$
\begin{subequations} \label{eq:dcvfbfs}
  \begin{align}
    \tilde a(w,  v) & = \int_0^1 w v \ d x + \sum_{\ell = 1}^\alpha \eta_{b,\ell} h^{6\ell-1} \int_0^1 w^{(2\ell)} v^{(2\ell)} \ d x, \label{eq:dcvfbfsa} \\
    \tilde b(w,  v) & = \int_0^1 w' v' \ d x + \sum_{\ell = 1}^\alpha \eta_{a,\ell} \pi^2 h^{6\ell-3} \int_0^1 w^{(2\ell)} v^{(2\ell)} \ d x, \label{eq:dcvfbfsb}
    % \tilde \ell(w) & = \ell(w),
    % + \sum_{\ell = 1}^\alpha \eta_{a,\ell} h^{6\ell-1} (\Delta^{\ell-1}  f, \Delta^{\ell} w)_{\partial \Omega}, \label{eq:dcvfbfsc} 
  \end{align}
\end{subequations}
and $\eta_{a,l}, \eta_{b,l}$ are penalty parameters set $\eta_{a,l} = \eta_{b,l} = 1$ by default.  With the one-dimensional bilinear forms defined above, using the tensor-product structure,  we define the 2D  bilinear forms as
\begin{equation} \label{eq:bp2d}
  \begin{aligned}
    \tilde{a}(\cdot, \cdot) & = \tilde{a}_x(\cdot, \cdot) \cdot \tilde{a}_y(\cdot,  \cdot), \\
    \tilde{b}(\cdot, \cdot)  & =  \tilde{b}_x(\cdot, \cdot) \cdot \tilde{a}_y(\cdot,  \cdot)+ \tilde{a}_x(\cdot, \cdot ) \cdot \tilde{b}_y(\cdot,  \cdot)
  \end{aligned}
\end{equation}
and the 3D bilinear forms as
\begin{equation} \label{eq:bp3d}
  \begin{aligned}
    \tilde{a}(\cdot, \cdot) & = \tilde{a}_x(\cdot, \cdot) \cdot \tilde{a}_y(\cdot,  \cdot) \cdot \tilde{a}_z(\cdot,  \cdot), \\
    \tilde{b}(\cdot, \cdot)  & =  \tilde{b}_x(\cdot, \cdot) \cdot \tilde{a}_y(\cdot,  \cdot)\cdot \tilde{a}_z(\cdot,  \cdot) + \tilde{a}_x(\cdot, \cdot ) \cdot \tilde{b}_y(\cdot,  \cdot) \cdot \tilde{a}_z(\cdot,  \cdot) + \tilde{a}_x(\cdot, \cdot ) \cdot \tilde{a}_y(\cdot,  \cdot) \cdot \tilde{b}_z(\cdot,  \cdot),
  \end{aligned}
\end{equation}
where $\tilde{a}_\xi, \tilde{b}_\xi, \xi = x,y,z$ are the bilinear forms in each dimension that we define similarly to~\eqref{eq:dcvfbfsa} and~\eqref{eq:dcvfbfsb}.  These bilinear forms lead to a new matrix problem
\begin{equation} \label{eq:mpfnew}
  \tilde M \ddot{\tilde U} + \tilde K \tilde U = F,
\end{equation}
which we solve using the explicit generalized-$\alpha$ method, see Section~\ref{sec:ga}.

\begin{remark}[Coercivity and optimality of~\eqref{eq:vfhnew}]
The boundary penalization bilinear terms are positive semi-definite. Thus, the new bilinear forms are coercive from the finite element theory, ensuring the stability of the spatial discretization. The discretization~\eqref{eq:vfhnew} is a Galerkin finite element method. The collection of B-spline basis functions forms a space that is a subspace of the usual finite element space. The a priori error estimate established for Galerkin finite elements also holds for the discretization~\eqref{eq:vfhnew}. As a consequence, we expect optimal convergence rates for the overall method to approximate a smooth solution $u$ to problem~\eqref{eq:pde}. That is,
\begin{align} \label{eq:ee}
  \| u(\cdot, T) - u^h(\cdot, T) \|_{L^2(\Omega)} &\le C(h^{p+1}+\tau^2), \\ 
  | u(\cdot, T) - u^h(\cdot, T)  |_{H^1(\Omega)} &\le C(h^p + \tau^2), 
\end{align}
where $C$ is independent of the mesh size $h$ and time-marching step size $\tau$.  We validate these estimates numerically in the next Section~\ref{sec:num}.
\end{remark}

\section{Numerical examples} \label{sec:num}

Numerical experiments validate the a priori error estimates~\eqref{eq:ee} for both the spatial and temporal discretization.  We then demonstrate the critical time step size increases in various examples.  The spectrum outliers appear when using isogeometric discretization with high-order elements; thus, we focus on $p\ge 3$ in our numerical experiments.

\begin{figure}[h!]
\centering
\includegraphics[height=6cm]{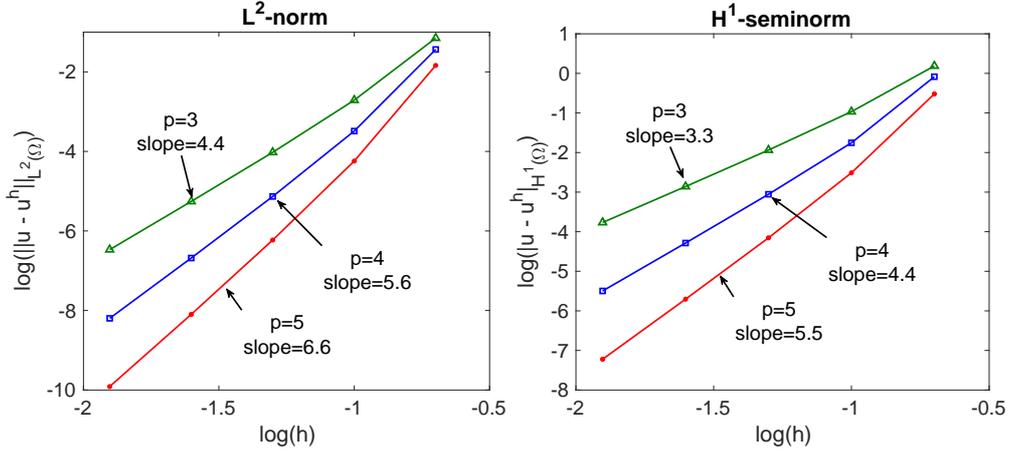} 
\caption{$L^2$-norm (left plot) and $H^1$-seminorm (right plot) errors of the boundary-penalized isogeometric analysis with $C^2$ cubic ($p=3$), $C^3$ quartic ($p=4$), and $C^4$ quintic ($p=5$) elements for the problem~\eqref{eq:pde} with $\kappa=1$ in 1D. }
\label{fig:errh1d}
\end{figure}

\subsection{Optimal approximation errors}

\begin{figure}[h!]
\includegraphics[height=6cm]{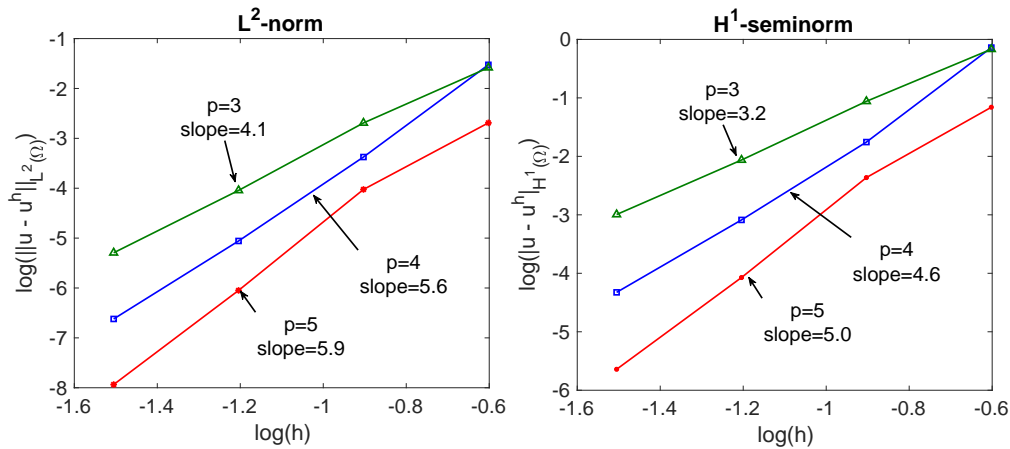} 
\caption{$L^2$-norm (left plot) and $H^1$-seminorm (right plot) errors of the boundary-penalized isogeometric analysis with $C^2$ cubic ($p=3$), $C^3$ quartic ($p=4$), and $C^4$ quintic ($p=5$) elements for the problem~\eqref{eq:pde} with $\kappa=1$ in 2D. }
\label{fig:errh2d}
\end{figure}

\subsubsection{Spatial approximation dominates the error}

 We first consider the test problem~\eqref{eq:pde} with $\kappa=1$ and a manufactured solution $u(x,t) = e^{t} \sin(3\pi x)$ in 1D and $u(x,t) = e^{t} \sin(3\pi x) \sin(3\pi y)$ in 2D. where we derive the forcing functions from~\eqref{eq:pde}.  We show how the spatial discretization may dominate the errors in a 1D problem; we set the final time $T=1$ with 10,000 time steps for the explicit generalized-$\alpha$ with $\rho=1$.  This setting guarantees that the spatial discretization dominates the error.  Figure~\ref{fig:errh1d} shows the $L^2$-norm and $H^1$-seminorm errors of the boundary-penalized isogeometric analysis for $p=\{3,4,5\}$.  The mesh sizes are $N = 5, 10, 20, 40, 80$.  We observe optimal error convergence rates which confirm the theoretical estimates in~\eqref{eq:ee} for the spatial discretization error dominance.

In 2D, the problem size of the spatial discretization increases quadratically with the discretization size in one dimension; we thus set the final time $T=0.01$ with 100 time steps to show the spatial accuracy. Figure~\ref{fig:errh2d} shows the $L^2$-norm and $H^1$-seminorm errors of the problem in 2D with mesh sizes $N = 4\times4, 8\times8, 16\times16, 32\times32$; the figure shows optimal spatial convergence.

\subsubsection{Temporal approximation dominates the error}

\begin{figure}[h!]
\centering
\includegraphics[height=6cm]{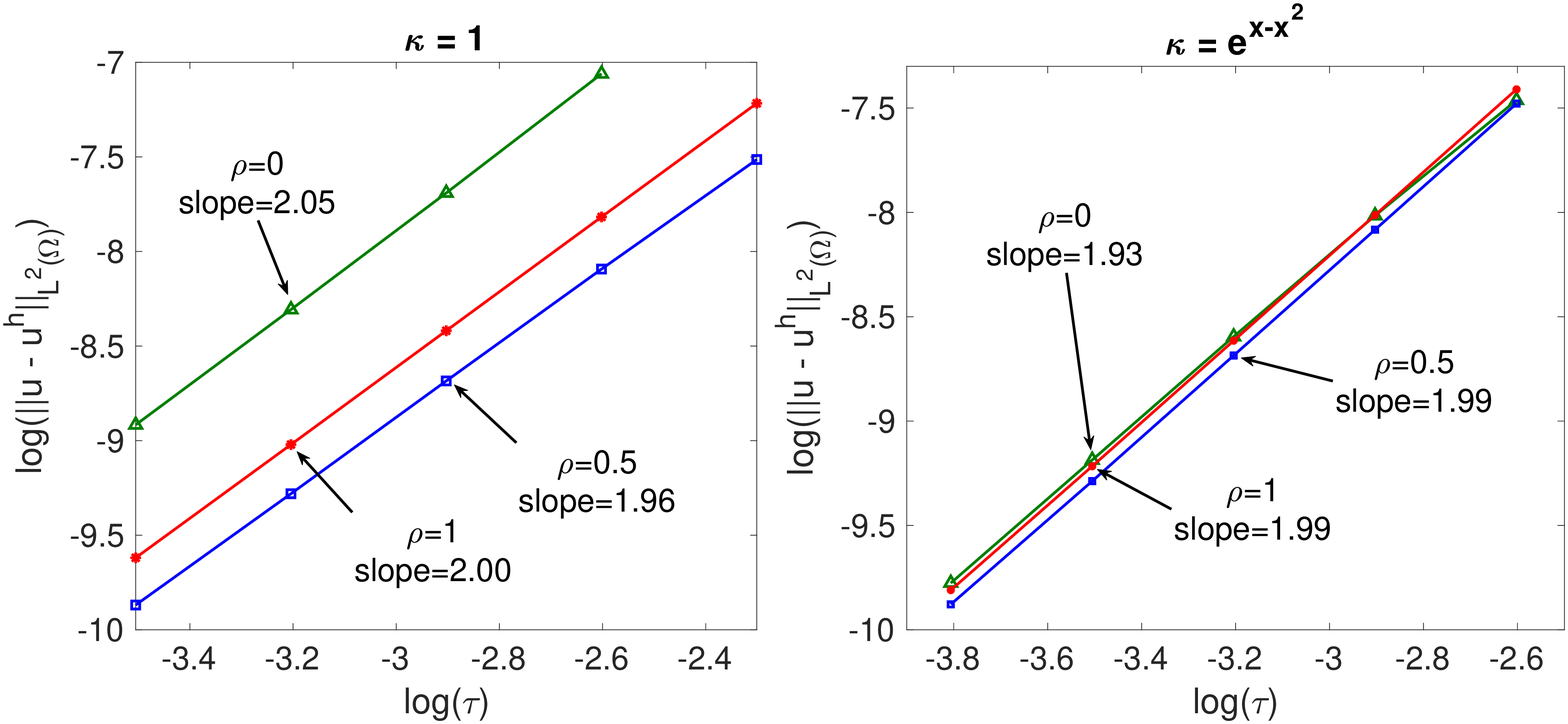} 
\caption{$L^2$-norm errors of the explicit generalized-$\alpha$ method for the problem~\eqref{eq:pde} with $\kappa=1$ (left plot) and $\kappa=e^{x-x^2}$(right plot) in 1D. }
\label{fig:errt1d}
\end{figure}

Now we focus on the $L^2$-norm error and analyze the limit where the temporal errors dominate. For this purpose, we apply the boundary-penalized isogeometric analysis with a fine mesh with $N=100$ in 1D. $C^4$ quintic ($p=5$) elements deliver a highly accurate spatial discretization so that the temporal discretization dominates the errors. Figure~\ref{fig:errt1d} shows the $L^2$-norm error for the explicit generalized-$\alpha$ method with $\rho=0,0.5,1$.  The left plot shows the errors for $\kappa=1$, while the right plot shows the errors for $\kappa=e^{x-x^2}$. We observe second-order accuracy in time which confirms the theoretical estimate~\eqref{eq:ee} for the temporal discretization by the explicit generalized-$\alpha$ method.

\subsection{Growth of the critical time-step size}

We consider the impact of boundary penalization on the critical time step sizes for the explicit generalized-$\alpha$ method.  In general, the critical time-step size depends on the parameter $\rho$ and the spatially discretized system's largest discrete eigenvalue (or frequency).  The conditional stability for the explicit generalized-$\alpha$ method is of the form
\begin{equation*}
  \tau \omega^h \le C_\rho,
\end{equation*}
being $C_\rho >0$ a constant independent of time-step $\tau$ and $\omega^h$ but dependent on the parameter $\rho$ of the scheme.  Then, the critical time step size becomes
\begin{equation*}
  \tau_c =  C_\rho/\omega^h_{\text{max}},
\end{equation*}
where $\omega^h_{\text{max}}$ is the largest frequency, that is, $\lambda^h_L = (\omega^h_{\text{max}})^2$ is the largest eigenvalue to~\eqref{eq:mevp}.  We perform a numerical study of critical time-step size using the examples from the subsections above.

\begin{figure}[h!]
\hspace{-0.5cm}
\includegraphics[height=6cm]{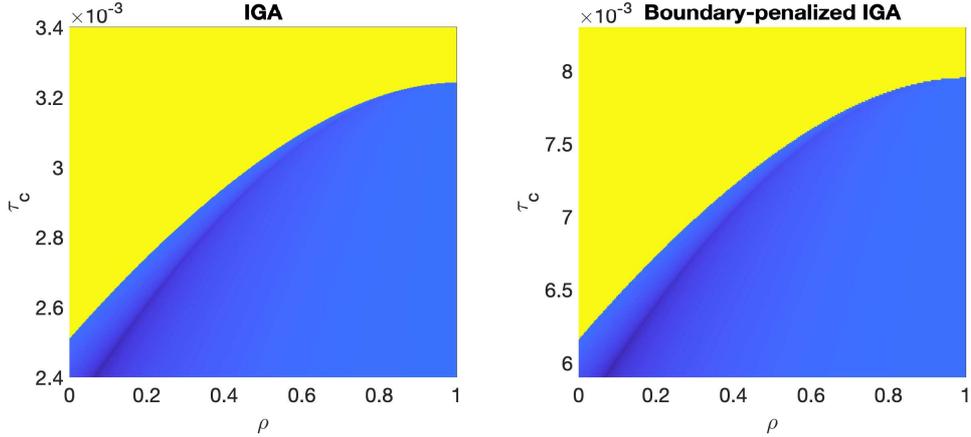} 
\caption{Stability regions with respect to the critical time step size $\tau_c$ and the parameter $\rho$. Left plot: region when using isogeometric analysis; right plot: region when using the boundary-penalized isogeometric analysis. The spatial discretization uses $C^5$ sextic B-spines with 80 uniform elements.}
\label{fig:stab}
\end{figure}

Table~\ref{tab:1d} compares the critical time steps for the explicit generalized-$\alpha$ method with $\rho=0$ for~\eqref{eq:pde} with $\kappa=1$ in 1D.  There are $N=5,10,20,40,80$ isogeometric elements with $p=\{3,4,5,6\}$ and maximal continuity.  The parameters with tildes denote the corresponding values when using boundary-penalized isogeometric elements.  The critical time-step size grows by a factor of about $\rho_{c,p} = \sqrt{\frac{(p-1)(p-2)}{4} + 1}$ for $p=\{3,4,5,6\}$.  For $p=1,2$, this formula also holds true as the boundary-penalized isogeometric analysis~\eqref{eq:vfhnew} reduces to the standard isogeometric analysis~\eqref{eq:vfh} that leads to the same critical time step size.  Figure~\ref{fig:stab} shows the comparison of the stability regions when using both the standard and  boundary-penalized isogeometric analysis.  Herein, we utilize $C^5$ sextic B-spline basis functions with 80 elements.  We observe that the critical time step size increases nonlinearly as the parameter $\rho$ increases.  For the largest critical time step size, one sets $\rho=1$.  For all values of $\rho$, the critical time-step size increases when using the boundary penalization technique to improve the performance of isogeometric analysis.

\begin{table}[h!]
\centering 
\begin{tabular}{ | c | c | cc | cc | c |}
\hline
$p$ & $N$ & \multicolumn{2}{c|}{$\lambda^h_L$ \qquad\quad $\tilde \lambda^h_L$}  & \multicolumn{2}{c|}{$\tau_c$ \qquad\quad $\tilde\tau_c$} & $\tilde \tau_c / \tau_c$  \\[0.1cm] \hline	
	&5	&402.8	&246.9	&7.72E-2	&9.86E-2	&1.28 \\
	&10    &1473.6	&987.5	&4.04E-2	&4.93E-2	&1.22 \\
3	&20	&5823.5	&3950.1	&2.03E-2	&2.46E-2	&1.21 \\
	&40	&23289.6	&15800.4	&1.02E-2	&1.23E-2	&1.21 \\
	&80	&93158.2	&63202.2	&5.08E-3	&6.16E-3	&1.21 \\ \hline
	&5	&680.9	&246.8	&5.94E-2	&9.86E-2	&1.66 \\
	&10	&2473.6	&987.2	&3.11E-2	&4.93E-2	&1.58 \\
4	&20	&9797.3	&3948.6	&1.57E-2	&2.47E-2	&1.58 \\
	&40	&39184.6	&15794.5	&7.83E-3	&1.23E-2	&1.58 \\
	&80	&156738.5	&63177.9	&3.91E-3	&6.16E-3	&1.58 \\ \hline
	&5	&1105.5	&246.8	&4.66E-2	&9.86E-2	&2.12 \\
	&10	&3976.8	&987.5	&2.46E-2	&4.93E-2	&2.01 \\
5	&20	&15722.0	&3952.3	&1.24E-2	&2.46E-2	&1.99 \\
	&40	&62874.0	&15845.2	&6.18E-3	&1.23E-2	&1.99 \\
	&80	&251495.8	&63894.6	&3.09E-3	&6.13E-3	&1.98 \\ \hline
	&5	&1703.9	&246.8	&3.75E-2	&9.86E-2	&2.63 \\
	&10	&6040.7	&987.2	&1.99E-2	&4.93E-2	&2.47 \\
6	&20	&23810.0	&3949.2	&1.00E-2	&2.47E-2	&2.46 \\
	&40	&95199.4	&15802.2	&5.02E-3	&1.23E-2	&2.45 \\
	&80	&380797.4	&63297.0	&2.51E-3	&6.16E-3	&2.45 \\ \hline
\end{tabular}
\caption{Comparison of the critical time steps for the explicit generalized-$\alpha$ method with $\rho=0$ when using isogeometric and boundary-penalized isogeometric discretizations for~\eqref{eq:pde} with $\kappa=1$ in 1D.}
\label{tab:1d} 
\end{table}

Lastly, Table~\ref{tab:1dk} shows the critical time-step sizes for~\eqref{eq:pde} with $\kappa=e^{x-x^2}$ in 1D while Table~\ref{tab:2d} shows the critical time step sizes for~\eqref{eq:pde} with $\kappa=1$ in 2D.  For the case with non-constant diffusion coefficient, the critical step size increment factor
$$\rho_{c,p}=\sqrt{\frac{(p-1)(p-2)}{4} + 1}$$
remains valid for coarse grids. For finer meshes, this factor decreases and tends to
$$\rho_{c,p}=p-2$$
for $p=\{3,4,5,6\}.$ For the 2D cases, the critical time-step size
increases by a similar factor, that is,
$$\rho_{c,p}=\sqrt{\frac{(p-1)(p-2)}{4} + 1}.$$
In summary, the boundary penalization technique increases the critical time-step size of the explicit generalized-$\alpha$ method in the isogeometric analysis of the second-order hyperbolic equation.

\begin{table}[h!]
\centering 
\begin{tabular}{ | c | c | cc | cc | c |}
\hline
$p$ & $N$ & \multicolumn{2}{c|}{$\lambda^h_L$ \qquad\quad $\tilde\lambda^h_L$}  & \multicolumn{2}{c|}{$\tau_c$ \qquad\quad $\tilde\tau_c$} & $\tilde \tau_c / \tau_c$  \\[0.1cm] \hline
	&5	&449.9	&247.0	&8.80E-2	&1.19E-1	&1.35 \\
	&10	&1588.6	&1071.7	&4.68E-2	&5.70E-2	&1.22 \\
3	&20	&6058.7	&4753.3	&2.40E-2	&2.71E-2	&1.13 \\
	&40	&23762.8	&19687.5	&1.21E-2	&1.33E-2	&1.10 \\
	&80	&94105.8	&79970.9	&6.08E-3	&6.60E-3	&1.08 \\ \hline
	&5	&732.5	&292.3	&6.90E-2	&1.09E-1	&1.58 \\
	&10	&2578.1	&1176.4	&3.68E-2	&5.44E-2	&1.48 \\
4	&20	&10009.4	&4778.8	&1.87E-2	&2.70E-2	&1.45 \\
	&40	&39616.1	&19596.4	&9.38E-3	&1.33E-2	&1.42 \\
	&80	&157607.6	&79726.3	&4.70E-3	&6.61E-3	&1.41 \\ \hline
	&5	&1164.8	&246.9	&5.47E-2	&1.19E-1	&2.17 \\
	&10	&4091.2	&1020.5	&2.92E-2	&5.84E-2	&2.00 \\
5	&20	&15954.0	&4653.7	&1.48E-2	&2.74E-2	&1.85 \\
	&40	&63345.2	&19515.3	&7.42E-3	&1.34E-2	&1.80 \\
	&80	&252445.0	&79579.1	&3.71E-3	&6.62E-3	&1.78 \\ \hline
	&5	&1773.4	&292.2	&4.43E-2	&1.09E-1	&2.46 \\
	&10	&6171.1	&1174.8	&2.38E-2	&5.45E-2	&2.29 \\
6	&20	&24073.3	&4755.8	&1.20E-2	&2.71E-2	&2.25 \\
	&40	&95733.3	&19478.2	&6.03E-3	&1.34E-2	&2.22 \\
	&80	&381872.3	&79452.8	&3.02E-3	&6.62E-3	&2.19 \\ \hline
\end{tabular}
\caption{Comparison of the critical time steps for the explicit generalized-$\alpha$ method with $\rho=0.5$ when using isogeometric and boundary-penalized isogeometric discretizations for~\eqref{eq:pde} with $\kappa=e^{x-x^2}$ in 1D.}
\label{tab:1dk} 
\end{table}

\begin{table}[h!]
\centering 
\begin{tabular}{ | c | c | cc | cc | c |}
\hline
$p$ & $N$ & \multicolumn{2}{c|}{$\lambda^h_L$ \qquad\quad $\tilde\lambda^h_L$}  & \multicolumn{2}{c|}{$\tau_c$ \qquad\quad $\tilde\tau_c$} & $\tilde \tau_c / \tau_c$  \\[0.1cm] \hline
	&5	&805.5	&493.8	&4.65E-2	&5.94E-2	&1.28 \\
	&10	&2947.3	&1975.0	&2.43E-2	&2.97E-2	&1.22 \\
3	&20	&11646.9	&7900.2	&1.22E-2	&1.48E-2	&1.21 \\
	&40	&46579.1	&31600.8	&6.12E-3	&7.42E-3	&1.21 \\
	&80	&186316.4	&126404.3	&3.06E-3	&3.71E-3	&1.21 \\ \hline
	&5	&1361.8	&493.6	&3.58E-2	&5.94E-2	&1.66 \\
	&10	&4947.2	&1974.3	&1.88E-2	&2.97E-2	&1.58 \\
4	&20	&19594.7	&7897.3	&9.43E-3	&1.49E-2	&1.58 \\
	&40	&78369.2	&31589.0	&4.71E-3	&7.43E-3	&1.58 \\
	&80	&313477.0	&126355.8	&2.36E-3	&3.71E-3	&1.58 \\ \hline
	&5	&2211.0	&493.7	&2.81E-2	&5.94E-2	&2.12 \\
	&10	&7953.6	&1975.0	&1.48E-2	&2.97E-2	&2.01 \\
5	&20	&31444.0	&7904.7	&7.44E-3	&1.48E-2	&1.99 \\
	&40	&125747.9	&31690.5	&3.72E-3	&7.41E-3	&1.99 \\
	&80	&502991.7	&127789.2	&1.86E-3	&3.69E-3	&1.98 \\ \hline
	&5	&3407.8	&493.6	&2.26E-2	&5.94E-2	&2.63 \\
	&10	&12081.3	&1974.4	&1.20E-2	&2.97E-2	&2.47 \\
6	&20	&47620.1	&7898.4	&6.05E-3	&1.49E-2	&2.46 \\
	&40	&190398.8	&31604.3	&3.02E-3	&7.42E-3	&2.45 \\
	&80	&761594.7	&126594.0	&1.51E-3	&3.71E-3	&2.45 \\ \hline	
\end{tabular}
\caption{Comparison of the critical time steps for the explicit generalized-$\alpha$ method with $\rho=0.5$ when using isogeometric and boundary-penalized isogeometric discretizations for~\eqref{eq:pde} with $\kappa=1$ in 2D.}
\label{tab:2d} 
\end{table}

\section{Concluding remarks} \label{sec:conclusion}

We study the impact of a boundary penalization technique for outlier removal on the resulting critical time-step size for isogeometric analysis with the explicit generalized-$\alpha$ method.  The boundary penalization technique reduces unphysical overshooting of the largest eigenvalues of the spatial discretization; consequently, the penalization increases the critical time-step sizes.  For the second-order hyperbolic equation~\eqref{eq:pde}, the spatially discrete eigenvalue $\lambda^h$ and the time step size $\tau$ are grouped into a form of $\lambda^h \tau^2$.  The boundary penalization technique reduces the largest eigenvalue by a factor of $\frac{(p-1)(p-2)}{4} + 1$, which then leads to an increment on the critical time step size by a factor of $\sqrt{\frac{(p-1)(p-2)}{4} + 1}$.  In the case of parabolic equation such as the heat equation $\dot u - \Delta u = f$, the discrete eigenvalue $\lambda^h$ and the time step size $\tau$ are grouped into a form of $\lambda^h \tau$, we thus expect an increment on the critical time step size by a factor of $\frac{(p-1)(p-2)}{4} + 1$.

With this in mind, one direction for future work would be to study the performance of the boundary penalization technique for solving the high-order initial-value second-order boundary-value problems.  For a partial differential equation with $k$-th order in time and second-order in space, one expects a critical time step size increment factor of $\Big( \frac{(p-1)(p-2)}{4} + 1 \Big)^{k/2}$.  An application of the proposed method is to use it to speed up the inverse modeling problems where a forward model is solved many times.

\section*{Acknowledgments} 

This publication was made possible in part by the Professorial Chair in Computational Geoscience at Curtin University and the Mineral Resources Business Unit of the Commonwealth Scientific Industrial Research Organisation, CSIRO, of Australia. This project has received funding from the European Union's Horizon 2020 research and innovation program under the Marie Sklodowska-Curie grant agreement No 777778 (MATHROCKS). The Curtin Corrosion Centre and the Curtin Institute for Computation kindly provide ongoing support.

%\section*{References}

%\bibliographystyle{plain}
%\bibliographystyle{siam}
%\bibliography{ref.bib}
%\bibliographystyle{siamplain}
\bibliographystyle{siam}
\bibliography{igaref}

\begin{thebibliography}{10}

\bibitem{behnoudfar2022explicit}
{\sc P.~Behnoudfar, G.~Loli, A.~Reali, G.~Sangalli, and V.~M. Calo}, {\em
  Explicit high-order generalized-$\alpha$ methods for isogeometric analysis of
  structural dynamics}, Computer Methods in Applied Mechanics and Engineering,
  389 (2022), p.~114344.

\bibitem{buffa2010isogeometric}
{\sc A.~Buffa, C.~De~Falco, and G.~Sangalli}, {\em Isogeometric analysis: new
  stable elements for the {S}tokes equation}, International Journal for
  Numerical Methods in Fluids,  (2010).

\bibitem{butcher2016numerical}
{\sc J.~C. Butcher}, {\em Numerical methods for ordinary differential
  equations}, John Wiley \& Sons, 2016.

\bibitem{calo2019dispersion}
{\sc V.~Calo, Q.~Deng, and V.~Puzyrev}, {\em Dispersion optimized quadratures
  for isogeometric analysis}, Journal of Computational and Applied Mathematics,
  355 (2019), pp.~283--300.

\bibitem{chung1993time}
{\sc J.~Chung and G.~Hulbert}, {\em A time integration algorithm for structural
  dynamics with improved numerical dissipation: the generalized-$\alpha$
  method}, Journal of Applied Mechanics, 60 (1993), pp.~371--375.

\bibitem{cottrell2009isogeometric}
{\sc J.~A. Cottrell, T.~J.~R. Hughes, and Y.~Bazilevs}, {\em Isogeometric
  analysis: toward integration of {CAD} and {FEA}}, John Wiley \& Sons, 2009.

\bibitem{cottrell2006isogeometric}
{\sc J.~A. Cottrell, A.~Reali, Y.~Bazilevs, and T.~J.~R. Hughes}, {\em
  Isogeometric analysis of structural vibrations}, Computer methods in applied
  mechanics and engineering, 195 (2006), pp.~5257--5296.

\bibitem{courant2008methods}
{\sc R.~Courant and D.~Hilbert}, {\em Methods of mathematical physics: partial
  differential equations}, John Wiley \& Sons, 2008.

\bibitem{craig2006fundamentals}
{\sc R.~R. Craig~Jr and A.~J. Kurdila}, {\em Fundamentals of structural
  dynamics}, John Wiley \& Sons, 2006.

\bibitem{Dalcin:2016}
{\sc L.~Dalcin, N.~Collier, P.~Vignal, A.~C{\^o}rtes, and V.~M. Calo}, {\em
  Pet{IGA}: A framework for high-performance isogeometric analysis}, Computer
  Methods in Applied Mechanics and Engineering, 308 (2016), pp.~151--181.

\bibitem{de1978practical}
{\sc C.~De~Boor}, {\em A practical guide to splines}, vol.~27, Springer-Verlag
  New York, 1978.

\bibitem{deng2022isogeometric}
{\sc Q.~Deng}, {\em Isogeometric analysis of bound states of a quantum
  three-body problem in {1D}}, in International Conference on Computational
  Science, Springer, accepted, 2022.

\bibitem{deng2018dispersion}
{\sc Q.~Deng, M.~Barto{\v{n}}, V.~Puzyrev, and V.~Calo}, {\em
  Dispersion-minimizing quadrature rules for {C}1 quadratic isogeometric
  analysis}, Computer Methods in Applied Mechanics and Engineering, 328 (2018),
  pp.~554--564.

\bibitem{deng2018ddm}
{\sc Q.~Deng and V.~Calo}, {\em Dispersion-minimized mass for isogeometric
  analysis}, Computer Methods in Applied Mechanics and Engineering, 341 (2018),
  pp.~71--92.

\bibitem{deng2021boundary}
{\sc Q.~Deng and V.~M. Calo}, {\em A boundary penalization technique to remove
  outliers from isogeometric analysis on tensor-product meshes}, Computer
  Methods in Applied Mechanics and Engineering, 383 (2021), p.~113907.

\bibitem{deng2021outlier}
\leavevmode\vrule height 2pt depth -1.6pt width 23pt, {\em Outlier removal for
  isogeometric spectral approximation with the optimally-blended quadratures},
  in International Conference on Computational Science, Springer, 2021,
  pp.~315--328.

\bibitem{deng2018isogeometric}
{\sc Q.~Deng, V.~Puzyrev, and V.~Calo}, {\em Isogeometric spectral
  approximation for elliptic differential operators}, Journal of Computational
  Science,  (2018).

\bibitem{deng2019optimal}
\leavevmode\vrule height 2pt depth -1.6pt width 23pt, {\em Optimal spectral
  approximation of 2n-order differential operators by mixed isogeometric
  analysis}, Computer Methods in Applied Mechanics and Engineering, 343 (2019),
  pp.~297--313.

\bibitem{evans2013isogeometric}
{\sc J.~A. Evans and T.~J. Hughes}, {\em Isogeometric divergence-conforming
  {B}-splines for the {D}arcy--{S}tokes--{B}rinkman equations}, Mathematical
  Models and Methods in Applied Sciences, 23 (2013), pp.~671--741.

\bibitem{evans2010partial}
{\sc L.~C. Evans}, {\em Partial differential equations}, vol.~19 of Graduate
  Studies in Mathematics, American Mathematical Society, Providence, RI,
  second~ed., 2010.

\bibitem{Garcia:2019}
{\sc D.~Garcia, D.~Pardo, and V.~M. Calo}, {\em Refined isogeometric analysis
  for fluid mechanics and electromagnetics}, Computer Methods in Applied
  Mechanics and Engineering, 356 (2019), pp.~598--628.

\bibitem{gomez2008isogeometric}
{\sc H.~G{\'o}mez, V.~M. Calo, Y.~Bazilevs, and T.~J.~R. Hughes}, {\em
  Isogeometric analysis of the {C}ahn--{H}illiard phase-field model}, Computer
  methods in applied mechanics and engineering, 197 (2008), pp.~4333--4352.

\bibitem{Hashemian:2021}
{\sc A.~Hashemian, D.~Pardo, and V.~M. Calo}, {\em Refined isogeometric
  analysis for generalized {H}ermitian eigenproblems}, Computer Methods in
  Applied Mechanics and Engineering, 381 (2021), p.~113823.

\bibitem{hiemstra2021removal}
{\sc R.~R. Hiemstra, T.~J. Hughes, A.~Reali, and D.~Schillinger}, {\em Removal
  of spurious outlier frequencies and modes from isogeometric discretizations
  of second-and fourth-order problems in one, two, and three dimensions},
  Computer Methods in Applied Mechanics and Engineering, 387 (2021), p.~114115.

\bibitem{hughes2012finite}
{\sc T.~J. Hughes}, {\em The finite element method: linear static and dynamic
  finite element analysis}, Courier Corporation, 2012.

\bibitem{hughes2005isogeometric}
{\sc T.~J.~R. Hughes, J.~A. Cottrell, and Y.~Bazilevs}, {\em Isogeometric
  analysis: {CAD}, finite elements, {NURBS}, exact geometry and mesh
  refinement}, Computer methods in applied mechanics and engineering, 194
  (2005), pp.~4135--4195.

\bibitem{hughes2014finite}
{\sc T.~J.~R. Hughes, J.~A. Evans, and A.~Reali}, {\em Finite element and
  {NURBS} approximations of eigenvalue, boundary-value, and initial-value
  problems}, Computer Methods in Applied Mechanics and Engineering, 272 (2014),
  pp.~290--320.

\bibitem{hughes2008duality}
{\sc T.~J.~R. Hughes, A.~Reali, and G.~Sangalli}, {\em Duality and unified
  analysis of discrete approximations in structural dynamics and wave
  propagation: comparison of p-method finite elements with k-method {NURBS}},
  Computer methods in applied mechanics and engineering, 197 (2008),
  pp.~4104--4124.

\bibitem{labanda2021explicit}
{\sc N.~A. Labanda, P.~Behnoudfar, and V.~M. Calo}, {\em An explicit
  predictor/multicorrector time marching with automatic adaptivity for
  finite-strain elastodynamics}, arXiv preprint arXiv:2111.07011,  (2021).

\bibitem{larsson2003partial}
{\sc S.~Larsson and V.~Thom{\'e}e}, {\em Partial differential equations with
  numerical methods}, vol.~45, Springer, 2003.

\bibitem{manni2022application}
{\sc C.~Manni, E.~Sande, and H.~Speleers}, {\em Application of optimal spline
  subspaces for the removal of spurious outliers in isogeometric
  discretizations}, Computer Methods in Applied Mechanics and Engineering, 389
  (2022), p.~114260.

\bibitem{newmark1959method}
{\sc N.~M. Newmark}, {\em A method of computation for structural dynamics},
  Journal of the engineering mechanics division, 85 (1959), pp.~67--94.

\bibitem{nguyen2015isogeometric}
{\sc V.~P. Nguyen, C.~Anitescu, S.~P. Bordas, and T.~Rabczuk}, {\em
  Isogeometric analysis: an overview and computer implementation aspects},
  Mathematics and Computers in Simulation, 117 (2015), pp.~89--116.

\bibitem{paz2012structural}
{\sc M.~Paz}, {\em Structural dynamics: theory and computation}, Springer
  Science \& Business Media, 2012.

\bibitem{piegl2012nurbs}
{\sc L.~Piegl and W.~Tiller}, {\em The {NURBS} book}, Springer Science \&
  Business Media, 1997.

\bibitem{puzyrev2017dispersion}
{\sc V.~Puzyrev, Q.~Deng, and V.~M. Calo}, {\em Dispersion-optimized quadrature
  rules for isogeometric analysis: modified inner products, their dispersion
  properties, and optimally blended schemes}, Computer Methods in Applied
  Mechanics and Engineering, 320 (2017), pp.~421--443.

\bibitem{Sarmiento:2017}
{\sc A.~F. Sarmiento, A.~M. C{\^o}rtes, D.~Garcia, L.~Dalcin, N.~Collier, and
  V.~M. Calo}, {\em Pet{IGA-MF}: a multi-field high-performance toolbox for
  structure-preserving {B}-splines spaces}, Journal of Computational Science,
  18 (2017), pp.~117--131.

\bibitem{serway2018physics}
{\sc R.~A. Serway and J.~W. Jewett}, {\em Physics for scientists and
  engineers}, Cengage learning, 2018.

\bibitem{tedesco2000structural}
{\sc J.~Tedesco, W.~G. McDougal, and C.~A. Ross}, {\em Structural dynamics},
  Pearson Education London, UK, 2000.

\bibitem{verhoosel2011isogeometric}
{\sc C.~V. Verhoosel, M.~A. Scott, T.~J. Hughes, and R.~De~Borst}, {\em An
  isogeometric analysis approach to gradient damage models}, International
  Journal for Numerical Methods in Engineering, 86 (2011), pp.~115--134.

\bibitem{wang2018structural}
{\sc Y.~Wang, Z.~Wang, Z.~Xia, and L.~H. Poh}, {\em Structural design
  optimization using isogeometric analysis: a comprehensive review}, Computer
  Modeling in Engineering \& Sciences, 117 (2018), pp.~455--507.

\bibitem{wanner1996solving}
{\sc G.~Wanner and E.~Hairer}, {\em Solving ordinary differential equations
  II}, vol.~375, Springer Berlin Heidelberg, 1996.

\bibitem{zhang2007patient}
{\sc Y.~Zhang, Y.~Bazilevs, S.~Goswami, C.~L. Bajaj, and T.~J. Hughes}, {\em
  Patient-specific vascular {NURBS} modeling for isogeometric analysis of blood
  flow}, Computer methods in applied mechanics and engineering, 196 (2007),
  pp.~2943--2959.

\end{thebibliography}

%\appendix{}
%\section{Other possibilities} 
%We present other 

\end{document}